\documentstyle[12pt,amssymb,amscd]{amsart} 
\def \part{{\bf part}}
\def\shpart{{\bf Shpart}}
\def\ess{{\bf ess}}
\def \Sf{{\bf S}}
\def\limdir{{\rm limdir}}
\def \A{{\cal A}}
\def \p{{\rm\bf p}}
\def\b {\bullet}
\def\grt{{\frak g}{\frak r}{\frak t}}
\def\g{{\frak g}}
\def \Im{{\rm Im}\;}
\def\d{{\frak d}}
\def \M{{\rm Mor}\;}
\def \C{{\cal C}}
\def \B{{\bf B}}
\def \BB{\B(U(\g(\cdot)))}
\def\u{{\flat}}
\def\disjoint{\sqcup}
\def\Id{{\rm Id}}
\def\Ob{{\rm Ob}\;}
\def\Mor{{\rm Mor}\;}
\def\Hom{{\rm Hom}\;}
\def\grt{{\frak grt}}
\def\sur{{\rm\bf sur}}
\def\morsur{{\rm Mor}_{\sur}}
\def\Shsur{{\rm\bf Shsur}}
\def \comm{{\rm\bf comm }}
\def\ass{{\rm\bf assoc}}
\def\lie{{\bf \rm lie}}
\def\hocomm{{\mbox{\em  ho}\comm}}
\def\hoass{{\rm\bf hoass}}
\def\pcd{{\rm \bf pcd}}

\def\defo{{\rm \bf def}}
\def \e{e}
\def \h{{\bf h}}
\def \CC{{\C_\bullet(\g(\cdot))}}
\def\Z{\Bbb Z}
\def \sh {{\bf\rm sh}}
\def\CD{{\bf CD}}
\def \PCD{{\bf PCD}}
\def \P{{\bf P}}
\def\CP{{\bf CP}}
\def\CPCD{{\bf CPCD}}
\def \PaCD{\PCD}
\def\PC{{\bf PC}}
\def\O{{\cal O}}
\def \S{{\cal S}}
\title{Action of the Grothendieck-Teichm\"uller\\ group
 On the operad of Gerstenhaber algebras}
\author{Dimitri Tamarkin}

\begin{document}

\newtheorem{Theorem}{THEOREM}[section]

\newtheorem{theorem}{THEOREM}[section]

 \newtheorem{definition}[Theorem]{DEFINITION}

 \newtheorem{Prop}[Theorem]{PROPOSITION}

 \newtheorem{Proposition}[Theorem]{PROPOSITION}

 \newtheorem{Lemma}[Theorem]{LEMMA}

 \newtheorem{Corollary}[Theorem]{COROLLARY}

\newtheorem{Claim}[Theorem]{CLAIM}

\maketitle
\section{introduction}
V. Drinfeld \cite{D}  and D. Bar-Nathan \cite{BN} introduced  an action of the (associated graded variant of) Grothendieck-Teichm\"uller group
on a certain algebraic structure which is called "Chord diagrams". This action can be transfered to an action on 
a free resolution of the operad of Gerstenhaber algebras. This action induces a map from  the Lie algebra $\grt$
of Grothendieck-Teichm\"uller group to  the differential 
graded Lie algebra of deformations of the operad of Gerstenhaber algebras.
We prove that this map is injective on the level of homology. In  other words, {\em any non-trivial element of the 
Grothendieck-Teichm\"uller group acts 
on a free resolution of the operad of Gerstenhaber algebras 
in a homotopically non-trivial way.} 

We start with definitions of the just mentioned objects.

\section{Notations}
$k$ stands for a fixed field of characteristic 0.  $[n]$ stands for the 
set $\{1,2,\ldots,n\}$.
\section{Chord Diagrams and the associated graded Grothendieck-Teichm\"uller group} 
 
\subsection{Chord diagrams}

 \subsection{Lie algebras $\g(T)$} Let $T$ be a finite set with at least 2 elements.
Consider a Lie algebra $\g(T)$ generated by the variables
$t_{ij}$, $i,j\in T$, $i\neq j$, $t_{ij}=t_{ji}$ satisfying the relations
$[t_{ij},t_{kl}]=0$ if $i,j,k,l$ are all different and  $[t_{ij},t_{ik}+t_{jk}]=0$,
whenever $i,j,k$ are all different. For a finite set $S$ with less than 2 elements set $\g(S)=0$.
For an integer $n\geq 0$ set $\g(n)=\g([n])$.
\subsubsection{Dilations}\label{Dilation} For $x\in k$ Let $D_x:\g(T)\to \g(T)$ be the automorphism such that
$D_x(t_{ij})=xt_{ij}$.

\subsection{Functoriality}
\subsubsection{Direct image}\label{direct} Let $f:S\to T$ be an {\em injective} map of finite sets. Define a Lie algebra homomorphism
$f_*:\g(S)\to \g(T)$ by the formula 
$$
f_*(t_{ij})=t_{f(i)f(j)}.
$$
It is straightforward to check that the relations are preserved.

\subsubsection{Inverse image}\label{inverse}
 Let $f:S\to T$ be a morphism of finite sets. Define a Lie algebra homomorphism
$f^*:\g(T)\to \g(S)$ by 
$$
f^*(t_{ij})=\sum\limits_{f(p)=i;f(q)=j}t_{pq}.
$$
Again, it is straightforward to check that the relations are preserved.
\subsubsection{A property} Let
$$
S\stackrel f\to T\stackrel g\to R
$$
be a sequence of finite sets and their maps such that $f$ is injective and the image of $gf$ has at most 1 element.
Then $\Im f_*$ and $\Im g^*$ are Lie subalgebras in $\g(T)$. 
\begin{Claim}\label{commutation}
$$
[\Im f_*,\Im g^*]=0.
$$  
\end{Claim}
\subsection{Operadic structure}
\subsubsection{}  For the definition of operad which is used here see
\ref{operads}. A slightly different  {\em language} is going to be used;
traditionally, one defines an operad 
in a symmetric monoidal category $C$
as a collection of $S_n$-modules
$\O(n)$. We instead prefer to speak about functors from the groupoid of finite sets $\S$ to $C$. This language has been now used by many authors. We hope that the reader will see the advantages of such a ' coordinate-free' (V. Hinich) approach.

 \subsubsection{} Category $L$ of Lie algebras  with the tensor product
$\g\otimes_L{\frak h}=\g\oplus\frak h$
is symmetric monoidal. The unit in $L$ is 0.
\subsubsection{} The assignment $X\mapsto \g(X)$ is a functor $\S\to L$, $f\in \Mor\S$ acting  via the direct image map $f_*$.   
\subsubsection{} Let $X,Y$ be finite sets and let $x\in X$. Let $Z=X-\{x\}\disjoint Y$. Define the map
$i:Y\to Z$ to be the natural inclusion and $p:Z\to X$ to be such that $p|_{X-\{x\}}$ is the natural inclusion
and $p(Y)=x$. 
Set $\circ_x^{XY}:\g(X)\oplus \g(Y)\to \g(Z)$ to  be $p^*\oplus i_*$. Claim \ref{commutation} implies
that $\circ_x^{XY}$ is a morphism of Lie algebras.

Let $e$ be a $1$-element set
and $\u_e:0\to \g(e)$ be the only existing morphism.
\subsubsection{} The collection $(\g(\cdot);\circ;\u_{(\cdot)})$ is an operad with unit.
\subsubsection{Chord diagrams} Let $A$ be the category of complete filtered associative algebras with unit. It is symmetric monoidal with respect 
to the completed tensor product of algebras. Let $U:L\to A$ be the functor of taking the completed universal enveloping
algebra. The isomorphism $U(\g\oplus {\frak h})\cong U(\g)\otimes U({\frak h})$ endows $U$ with a tensor structure.
Therefore, $U(\g(\cdot))$ is an operad of complete filtered
associative algebras. Denote it by $\CD$. 
This operad is almost what we need, but it does not satisfy us because it has no
automorphisms except $D_x$ (see \ref{Dilation}). 

\subsubsection{} Let $Cat$ be the category of small complete filtered $k$-linear categories (morphisms between any two
fixed objects form a  complete filtered vector space and the composition is a continuous bilinear map). It is also symmetric monoidal:
$\Ob(C\otimes D)=\Ob C\times \Ob D$; $Mor(x\times y, x'\times y')=\Mor(x,x')\otimes \Mor(y,y')$.

Any associative algebra $X$ can be viewed as a category with a unique object $x$ with  its endomorphisms being $X$.   
This defines a symmetric monoidal functor $C:A\to Cat$ and an operad structure on $C(U(\g(\cdot)))$.  

\subsubsection{} For a finite set $X$ let $p(X)$  be the set of all non-associative monomials in $X$ such that
any element of $X$ enters once. Other names for $p(X)$ are  binary trees with the set of terminal vertices $X$
or parenthesized permutations of $X$. 

We can define $p(X)$ inductively by assigning: 
\begin{enumerate}
\item[-]$P(\emptyset)=\{\u\}$, where $\u$ is a fixed element;
\item[-]$p(e)=e$ for any 1-element set$e$;
\item[-]
$$
p(X)=\disjoint_{X=X_1\disjoint X_2} p(X_1)\times p(X_2),
$$
where $X_{1,2}$ are non-empty.
\end{enumerate}
The $p(\cdot)$ form an operad of sets. It is an operad generated by one binary operation (denote it by $\cdot$) and by one 0-ary operation 
$\u$. The relations are 
$$
\circ_1^{[2],\emptyset}(\cdot,\u)=\circ_2^{[2],\emptyset}(\cdot,\u)={\rm Id}.
$$
Let $p'(X)$ be the category whose objects are $p(X)$  and there is a unique morphism between any 2 objects.
Let $\P(X)$ be the linear span of $p'(X)$. Then $\P(\cdot)$ is an operad in $Cat$.
\subsubsection{} Set $\PCD(X)=\CD(X)\otimes \P(X)$. The $\PCD(\cdot)$ is an operad in $Cat$. This is the object we are going to
study.
\subsubsection{} Note that the space of morphisms between any two fixed objects 
of $\CD(X)$ as well as $\P(X)$ has  a natural cocommutative coalgebra structure.
This defines  a cocommutative and coassociative map of operads $\PCD\to \PCD\otimes \PCD$.

In other words, let $a$ be a commutative $k$-algebra. Let $CHopf_a$ be the category of small categories such that space of morphisms between any two objects
is a complete filtered co-commutative $a$-coalgebra and the composition is  compatible with these coalgebra structures.  
Then $\PaCD$ is an operad in $CHopf_k$.

\subsection{Grothendieck-Teichm\"uller group}
\subsubsection{Change of the base ring}  let $a$ be a commutative $k$-algebra. Let $\PaCD\otimes a$ be an operad in
$CHopf_a$ such that the objects in $\PaCD\otimes a(X)$ are the same as in $\PCD(X)$ and the morphisms in
$\PaCD\otimes a $ are obtained from the ones in $\PCD$ by taking the completed tensor product with $a$.
\subsubsection{}
Denote by $GRT(a)$ the group of automorphisms of $\PCD\otimes a$ in the category of operads in $CHopf_a$
which are identical on $\Ob\;\PCD\otimes a(X)$ for any $X$. 
Let $F^iGRT(a)$ be the subgroup of such automorphisms $F$ that for any $x,y\in Ob\;\PCD\otimes a$, we have
$$
(Id-F)({\rm Mor}\;(x,y))\subset F^i ({\rm Mor}\;(x,y)).
$$
\subsubsection{} We have a natural map ${\frak D}:a^\times \to GRT(a)$ determined by the dilations $D_x$ from  \ref{Dilation}
\subsubsection{} It is well known that $GRT$ is a pro-algebraic group. Let $\grt$ be its complete filtered 
Lie algebra. It can be also defined as a Lie algebra of derivations of $\PCD$ as an operad in $CHopf_k$.
The map $\frak D$ defines a map of Lie algebras $\d:k\to \grt$. The action of $\d(1)$ on $\grt$ is diagonalizable
with nonnegative integer eigenvalues.
Let 
$$\grt^i=\{x\in \grt:\  [\d(1),x]=ix \}.
$$
Then $\grt^{< 0}=0$; $\grt^0=\d(k)$; $F^i\grt=\prod_{k\geq i}\grt^k$. 

In particular, this implies that $GRT$ is a pro-nilpotent extension of $G_m$.
\subsubsection{Notation}In the sequel we denote
$$
\grt_i=F^i\grt.
$$
\subsubsection{} let $[n]=\{1,2\,\ldots,n \}$. Let  $b=t_{12}\in {\rm Mor}\;((12),(21))$ in
the category $\PCD\otimes a([2])$; Let $\phi=1\in {\rm Mor}\;(((12)3),(1(23))$ in $\PCD\otimes a([3])$.
An element $x\in G(a)$ is completely determined by prescription of $$x(b)=\alpha(x)t_{12}\in \g(2)\otimes a
$$
 and 
$$
x(\phi)\in (U(\g(3))\otimes a)_{\mbox{grouplike}}.
$$
	All tensor products here are completed (this makes difference 
if $a$ is infinite-dimensional).
  The map $x\mapsto \alpha(x)$ determines (a unique possible)
splitting 
$$
k\stackrel \d\to \grt\stackrel \alpha\to k,
$$
whose kernel is $\grt_1$.
\subsubsection{} \label{phigrt}Any element  $u\in \grt_1$ (viewed as a derivation of $\PCD$) is
completely determined by $u(\phi)\in \g(3)$. We will also denote
$u(\phi)$ by $\phi(u)$.

It is known what is the image $\grt_1(\phi)\subset \g(3)$. It is described by the following two conditions.
\subsubsection{The first condition}\label{first condition} It stems from the following chain of morphisms in $ \PCD([3])$ (each of morphisms is given
by $1\in U(\g(3))$):
$$
(12)3\to (21)3\to 2(13)\to 2(31)\to (23)1\to 1(23)\to (12)3.
$$
The composition of these six morphisms is identity. Applying a derivation $u$ to this chain, we obtain:
$$
((213)-(231)-(123)).u(\phi)=0,
$$
 here $.$ means the action of the symmetric group. 
Similarly, one obtains:
$$
((132)-(312)-(123)).u(\phi)=0.
$$
This condition follows from the former one by acting on it with an appropriate element of the group algebra of $S_3$.

Altogether these conditions mean that {\em $u(\phi)$ vanishes on shuffles} (see \ref{shuffles}).
\subsubsection{Second condition} \label{second condition}It comes from the following chain of morphisms:
$$
((12)3)4\to (1(23))4\to 1((23)4)\to 1(2(34))\to (12)(34)\to ((12)3)4,
$$
each morphism is again 1.

To describe the corresponding condition introduce the following notation for surjective  nondecreasing morphisms 
$[4]\to [3]$: for such a morphism we will enclose into brackets the two elements that go to the same element:
$(12)34$ stands for a morphism that sends $1,2$ to $1$; $3,4$ to $2,3$ respectively.

Similarly, we describe an injective increasing morphism $[3]\to [4]$ by specifying its image:
$134$ sends 1,2,3 to 1,3,4 respectively.

Then the condition expresses itself as follows:
$$
\large\{(123)_*+(1(23)4)^*+(234)_*-(12(34))^*-((12)34)^*\large\}u(\phi)=0.
$$
\subsubsection{} Any element satisfying the above conditions is of the form $u(\phi)$ for a unique $u\in \grt_1$.

\subsubsection{}\label{quotient}
Let $c=1\in {\rm Mor}\;((12),(21))$ in $\PCD\otimes a([2])$. One can show that $c$ is fixed by any
element of $GRT(a)$. 
\subsection{ A modification} This subsection will not be used in the sequel.
\subsubsection{} Let $X$ be a finite set. 
Set $pc(X)$ to be the set of all commutative non-associative monomials in $X$
In other words, let $pc(\cdot)$ be the free operad generated by one  binary commutative operation and by one 0-ary operation representing the unit.
Let $\PC(\dot)$ be the corresponding operad in $CHopf_k$. One has the quotient map of operads
$\P(\dot)\to \PC(\dot)$ resulting in a  map $q: \PCD\to \CD\otimes \CP$. Set $\CPCD=\CD\otimes \CP$.
\subsubsection{} It follows from \ref{quotient} that the $GRT$ action on $\PCD$ naturally
descends to an action on $\CPCD$.
\subsection{Nerves}
\subsubsection{Definition} Let $C$ be a category from $Cat$. 
Let $\u$ be the category with one object and 1-dimensional space of morphisms. Assume that an augmentation $\varepsilon:C\to \u$ is given.

Define the simplicial complete filtered vector space $NC_\b$
such that 
\begin{eqnarray*}
&NC_k&\nonumber\\ 
&=&\hat{\bigoplus\limits_{X_i\in {\rm Ob}\; C} }\M(X_0,X_1)\otimes \M(X_1,X_2)\otimes \cdots\otimes \M(X_{k-1},X_{k}),\nonumber\\
\end{eqnarray*} 
where $\hat{\cdot}$ means completion. The face map $d_0$ (resp $d_k$) is given by applying $\varepsilon$ to the first (resp.
to the last) tensor  factor, the remaining factors being unchanged. All the other face maps $d_i$, $i=1,\ldots,k-1$,
are given by applying  the composition of the $i$-th and $(i-1)$-th factor. The degeneration maps $s_i$ is the insertion of the identity 
morphism ${\rm Id}\in \M(X_i,X_i)$.

$N$ is a symmetric monoidal functor from the category $Cat'$
of augmented objects in $Cat$ to the category 
of complete filtered simplicial vector spaces.

The operads $C(U(\g(\cdot)))$ and $\PCD$ can be viewed as operads in $Cat'$, their nerves thus being
operads of simplicial complete filtered vector spaces. 

The functor $Ch$ of taking the reduced chain complex of a simplicial
vector space has a symmetric structure via the Eilenberg-Zilber map. Let $\C=ChN$. 

\subsubsection{ Agreement}\label{agreement} The reduced chain complex will be always put  
in the negative degrees so that the differential would raise the degree by 1.
 
\subsubsection{} By applying the functor $\C$ we obtain the following operads of complexes of complete filtered vector spaces:
$$
 \B(U(\g(\cdot)))=\C(C(U(\g(\cdot))));\ \pcd=\C(\PCD).
$$

Note that the $\B(U(\g(X))$ is the reduced bar complex for an augmented complete filtered associative algebra $U(\g(\cdot))$.
The $\grt$ acts on $\pcd$; the action of $\grt_0$ splits the filtration. 
 The action of $\grt_0$ is also naturally defined on $\B(U(\g(\cdot)))$ (but not of the whole $\grt$).

The natural maps of operads $\P\to \u$ induce  a quasi-isomorphism
\begin{equation}\label{q}
q:\pcd\to \B(U(\g(\cdot))).
\end{equation}
\section{Gerstenhaber algebras}
\subsection{Definition}
A Gerstenhaber algebra structure on a complex $V$ is a prescription of
a commutative associative product 
$$
\cdot:S^2V\to V
$$ and a Lie bracket
$$
\{.\}:\Lambda^2(V[1])\to V[1]
$$
 making $V[1]$ a differential graded Lie algebra.
Note that this bracket can be also described as a map $S^2V\to V$ of degree $-1$.
The product and the bracket must satisfy the Leibnitz identity:
$$
\{a,bc\}=\{a,b\}c+(-1)^{(|a|-1)|b|}b\{a,c\}.
$$ 
Let $e_2$ be the operad of Gerstenhaber algebras. We denote by $m$ (resp $b$) the generators
representing the commutative product and the bracket.
\subsection{A map $\mu:e_2\to \B(U(\g(\cdot)))$}

Note that 
$$\B(U(\g(\cdot)))^{-l}=(U(\cdot)/k)^{\otimes l}.
$$
(see \ref{agreement} for the explanation of the degree $-l$).
Let 
$$
t=t_{12}\in U(\g(2))/k\subset \B(U(\g(2)))^{-1}
$$ 
and
$$c=1\in U(\g(2))/k\subset \B(U(\g(2)))^{-1}.
$$
\begin{Claim} The assignment $\mu(m)=c$; $\mu(b)=t$ extends to a morphism $\mu:e_2\to  \B(U(\g(\cdot)))$
which is a quasi-isomorphism of operads.
\end{Claim}
\subsubsection{Recapitulation}We have  the following diagram of quasi-isomorphisms of operads:
$$
\pcd\to  \B(U(\g(\cdot)))\leftarrow e_2.
$$
The $\grt$ acts on $\pcd$.
\subsubsection{} The homotopy theory implies that this defines (uniquely up to a homotopy equivalence)
an $L_\infty$-map $m$ from $\grt$ to the deformation complex $\defo(e_2)$
of $e_2$, the latter being the differential graded Lie algebra of derivations of a free resolution of $e_2$.
(more exactly, of a cofibrant one \cite{H}). To be certain, one can choose  the Koszul resolution of $e_2$.
Though $m$ is defined up to a homotopy, the induced map on cohomology $H^0(m): \grt\to H^0(\defo(\e_2))$ is uniquely
defined. Our goal is to prove that 
\begin{Theorem}
$H^0(m)$ {\em is injective}
\end{Theorem}
\section{Proof of the theorem}
\subsection{Homotopy-theoretical background }
\subsubsection{Relative derivations}\label{derivations}
 Let $f:A\to B$ be a map of operads of complexes
of  vector spaces. Let $A',B'$ be the same operads
but with zero differential. Let $a_n=k[h_n]/(h_n^2)$ be the graded commutative ring such that $|h_n|=n$.
Recall that $A',B'$ are functors from the groupoid of finite sets to graded vector spaces.
A  morphism of functors
 $g(\cdot):A'(\cdot)\to B'(\cdot)$  such that each  $g(X)$ is of grade $n$ is called an $f$-derivation if
the morphism $f+h_{-n}g:A'\to B'\otimes_k a_{-n}$ is a morphism of operads over $k$. Denote by 
$D(f)'_n$ the vector space of all $f$-derivations. The direct sum with respect to all $n$ is a graded space.    
Denote it by $D(f)'$. 
\subsubsection{} Let $I_A:A\to A$ be an identity morphism. Then  $D(I_A)'$ has a natural structure of differential
graded Lie algebra with respect to the usual commutator. The differential on $A$ determines an element
$d_A\in D(I_A)'$ such that $[d_A,d_A]=0$.  
\subsubsection{}
Define a differential on $D(f)'$ by the formula
$$
dx=d_Bx-(-1)^{|x|}xd_A.
$$
Denote by $D(f)$ the corresponding complex.
\subsubsection{} Let $f_i:A_i\to A_{i+1}$, $i=1,2,\ldots n$ be chain of morphisms of operads.
One has an obvious map of complexes
$$
(f_1,\ldots,f_n):\bigoplus\limits_k D(f_k)\to D(f_nf_{n-1}\cdots f_1).
$$
Denote by $f_{1}^*$ the restriction
of $(f_1f_2)$ onto $D(f_1)$ and by $f_{2*}$ the restriction onto $D(f_2)$.  
\subsubsection{} Let $A$ be an operad and $A_i\to A$ be cofibrant resolutions
\cite{H}.
Then the differential graded Lie algebras $D(I_{A_i})$ are   quasi-isomorphic, the quasi-isomorphisms
are uniquely defined up to a homotopy. 
Each of $D(I_{A_i})$ is called {\em the deformation Lie algebra of $A$}. 

Furthermore, there exist uniquely defined  up to a homotopy  $L_\infty$-maps 
$\gamma_i:D(I_A)\to D(I_{A_i})$ for each $i$. In particular, the cohomology 
$$
H^\bullet_i=H^\bullet(D(I_{A_i}))
$$
are canonically identified and there exists a canonic map
$H^\bullet(D_A)\to H^\bullet_i$.
\subsubsection{} Let $f:B\to A$ be a map of operads such that $B$ is cofibrant. We have a map
$f^*:D(I_A)\to D(f)$.  Let $A_1\to A$ be a cofibrant resolution of $A$. 
\begin{Claim} There exists an $L_\infty$-map $k: D(I_{A_1})\to D(f)$ such that the composition
$k\gamma_1:D(I_A)\to D(f)$ is homotopy equivalent to $f^*$.
\end{Claim}  
\begin{Corollary}
If $H^{\bullet}f^*$ is injective, then so is $H^{\bullet}\gamma_1$.
\end{Corollary}
Let $l:A\to C$ be another map. We have a map $l_*: D(f)\to D(lf)$,
the injectivity of $H^{\bullet}(l_*f^*)$ will then also imply the injectivity of $\gamma_i$ .
If $\alpha\subset D(I_A)$ is a subcomplex, then the injectivity of
 $H^\b(\gamma_j)$ on $H^\b\alpha$
implies the injectivity of $H^\b\gamma_1$ on $H^\b\alpha$.
  
\subsection{Back to GT} let $\p=\C(\P)$ be an operad of reduced chains of the nerve of the operad $\P$  
of parenthesized permutations. The natural map $\P\to \PCD$ induces a map
$i:\p\to\pcd$. Also we have a map $q:\pcd \to \BB$ (see Eq. (\ref{q})). The $\grt$-action on $\pcd$ induces a map
$\grt\to D(I_{\pcd})$. We have a map 
$$\pi:\grt\to D(I_\pcd)
\stackrel {i^*q_*}\to D(qi).
$$
\subsubsection{A property of $\pi$}
\begin{Claim} Let $u\in\grt$. Then $\pi(u)(\p^{s}(X))=0$ whenever $s<-1$.
\end{Claim}
 
\pf The space $\p^{s}(X)$ by definition has a basis each element of which
is specified by a sequence $\{x_0,\ldots, x_{-s}\}$ of objects of $\P(X)$ 
(recall that the objects of $\P(X)$  are parenthesized permutations)
such that for any $j$ we have
 $x_j\neq x_{j+1}$. of elements from $X$. This is so because the space of morphisms between any two objects is 1-dimensional and we assume
that the morphisms $x_i\to x_{i+1}$ are all 1's. The map $i$ will include these 1's into $U(g(X))$. The action
of a derivation can change at most one of these 1's, therefore, if $s<-1$, at least one of these 1 will remain, and
the application of $q$ will yield 0.
\endpf

Let 
$$
p'(X)=\tau_{\geq -1}p(X)=p(X)/(p^{\geq -1}(X)+dp(X)^{\geq -1})  
$$ 
be the quotient.
One checks that these quotients  inherit the operad structure from $p$ . 
 Let $s:p\to p'$ be the quotient map of operads.  Then the map $q$ passes through $s$:
$q=q's$.
We have a map $s^*:D(q')\to D(q)$
\begin{Corollary} The map $\pi$ passes through $s^*$: $\pi=s^*\pi'$ for a unique $\pi'$ 
\end{Corollary}
\subsubsection{} Let $\hocomm$ be the operad of homotopy commutative algebras (see \ref{hocomm}). There exists a quasi-isomorphism $c:\hocomm\to \p'$ such that $c(m_2)=((12)+(21))/2\in \p'([2])$.
The image of $m_3$ is uniquely determined by the conditions
$$
dc(m_3)+\frac 12\{c(m_2),c(m_2)\}=0.
$$
and
$c(m_l)=0$ for $l>3$. 
\subsubsection{} The map $c$ determines a map $c^*:\grt\to D_{q'c}$, where $q'c:\hocomm\to \BB$ is
the map which sends $m_2$ to $1\in U(\g_2)$ and all $m_l$, $l>2$ to zero.
\subsubsection{}\label{reduc} Our theorem now reduces to proof of injectivity of $c^*$.
\subsubsection{} We see that $c^*(u)(m_3)=\phi(u)\in\g([3])\subset U(\g([3]))/k$ (see \ref{phigrt}).
\subsection{Recapitulation} We will reformulate the reduction statement in \ref{reduc}
in  detail. The rest of the paper will be devoted to the proof of this statement.

Let $\BB$ be the operad of reduced bar complexes. Let $Q:\hocomm\to \BB$
be such that $Q(m_2)=1\in \BB^0$ and $Q(m_l)=0$ for all $l>2$.
Let $D(Q)$ be the space of all $Q$-derivations (see   \ref{derivations}).
Informally, $D(Q)$ is the set of all maps $\psi=\prod\psi_n$, where 
$\psi_n:\hocomm(n)\to \BB(n)$, satisfying the Leibnitz identity.
An element of $D(Q)$ is uniquely defined by its values on the generators 
$m_k\in \hocomm$. We have a map of complexes $f_*:\grt_1\to D(Q)$
such that $f_*(u)(m_k)=0$ for $k\neq 3$;
$f_*(u)(m_3)=\phi(u)$ (for $\phi$ see \ref{phigrt}).
We want to prove that $f_*$ induces an injective map on cohomology.
\section{Simplifications and reformulations}
\subsection{Reduction to an operad of chains}
\subsubsection{An operad of chains}
Let $(C_\bullet(\g(X)),\partial)$ be the Chevalley-Eilenberg chain complex of the Lie  
algebra $\g(X)$ .  The functor of taking chains is a symmetric monoidal functor
from the category of Lie algebras with tensor product being the direct sum of Lie algebras
to the category of complexes with the usual tensor product (actually to the category of
coalgebras but we do not need it). Therefore, $\CC$ is an operad of complexes
\subsubsection{}
We have an inclusion $i:C_\bullet({\g(X)})\to \B(\bar U(\g(X)))$ such that
$$
x_1\wedge\cdots \wedge x_m\mapsto \frac 1{m!} 
\sum\limits_{\sigma\in  S_m}x_{\sigma_1}\otimes\cdots\otimes x_{\sigma_m}
\in \B(\bar U(\g(X)))^{-1}.
$$
One sees that $i$ is a quasi-isomorphism of operads.
\subsubsection{} The map $Q:\hocomm\to \BB$ uniquely passes
through $\CC$: there exists a unique $R:\hocomm\to \CC$ such that
$Q=iR$.  We have an inclusion $i_*:D(R)\to D(Q)$. 
  The map
$c^*\grt_1\to D(Q)$ uniquely passes through $D(R)$ thus defining a map
$F:\grt_1\to D(R)$. Our theorem  is now reduced to this situation.
\subsubsection{Recapitulation} \label{reform}
The map $R:\hocomm\to \CC$ is defined by
$R(m_2)=1\in C_0(\g(2))$ and $R(m_k)=0$ for $k>2$.

Define $D(R)$ as explained in \ref{derivations}.
The map $F:\grt_1\to D(R)$ is defined by   $F(u)(m_l)=0$, 
 $l\neq 3$, $u\in\grt_1$. We set 
$$F(m_3)=\phi(u)\in \g(3)\subset
 C^{-1}(\g(3))\stackrel{\rm def}{=}C_1(\g(3)).
$$
We  want to prove that $F$ is injective on the level of cohomology.

\subsection{Another description of D(R)} In this section
we will see that to define $D(R)$ we do not need
to know an operad structure on $\CC$ but rather a weaker
structure of a {\em module} over the operad $\comm$. A similar 
statement holds in general situation too. 

\subsubsection{A module over an operad} Let $\O$ be an operad 
in a symmetric monoidal category $\C$
and let
$M$ be a functor from the groupoid   of finite sets $S$ to $\C$.
{\it A structure of an $\O$-module on $M$} is by definition a collection of composition maps
$$
\circ_x^{XY}:\O(X)\otimes M(Y)\oplus M(X)\otimes O(Y)
\to M(X-\{x\}\disjoint Y),
$$
where $X,Y$ are finite sets, $x\in X$.  These maps must be natural in $x,X,Y$ and 
satisfy the following property. Let 
\begin{eqnarray*}
(\circ_x^{XY})^e:\O(X)\otimes\O(Y)\oplus M(X)\otimes M(Y)\oplus\O(X)\otimes M(Y)\oplus M(X)\otimes O(Y)\\
\to M(X-\{x\}\disjoint Y)\oplus \O(X-\{x\}\disjoint Y)
\end{eqnarray*}
be the extended map such that

- on $\O(X)\otimes \O(Y)$ it is equal to the composition map in $\O$;

-on $M(X)\otimes M(Y)$ it is equal  to 0;

-on the remaining summands it is equal to $\circ_x^{XY}$.
Then $\circ^e$ must define an operadic law on $\O\oplus M$.
 \subsubsection{} Let $M_i$ be $\O_i$-modules. Then $\bigotimes_i M_i$
  is an $\bigotimes_i \O_i$-module 
\subsubsection{} Let $f:\O_1\to \O_2$ be a map of operads. It naturally defines
an $\O_1$-module structure on $\O_2$ in a way that
the corresponding map $\circ_x^{XY}=\circ_{\O_2,x}(f\otimes Id\oplus Id\otimes f)$.
In particular, any operad is a module over itself via the identity map.
\subsubsection{} Let $\hocomm\stackrel p\to \comm\stackrel g\to \O$
be a sequence of maps of operads, where $p$ is the natural projection.
 Our goal is to describe $D(gp)$ in terms of the $\comm$-module structure on
$\O$.  Note that the space $G(n)$ of $n$-ary  generators of $\hocomm$
 is identified with the dual space $\lie'(n)*$, where $\lie'=\lie\{-1\}$ is the shifted operad of Lie algebras.  Any derivation is uniquely described by
a map (of not necessarily 0-th grading) of functors from the groupoid of finite sets to complexes
$G\to \O$ or, by an element of 
$$
H=H(\O)=\bigoplus_{n=2}^{\infty} (\lie'(n)\otimes \O(n))^{S_n}.
$$
Let $\h(X)=\lie'(X)\otimes \O(X)$ be the tensor product of functors.
We have a $\lie'\otimes \comm$-structure on  $\h$. Note that
$\lie'\cong \lie'\otimes \comm$.
\subsubsection{} We have an inclusion of $S_n$ invariants
$H(n)\subset \h(n)$. Let $p_n:\h(n)\to H(n)$ be a unique equivariant splitting .  
\subsubsection{} Let $b\in\lie'(2)$ be the element representing the Lie bracket.
It has degree $1$. For an element $ x\in H(n)$  set 
$$
dx=p_{n+1}\large(\circ_1^{[n][2]}(x,b)-(-1)^{|x|}\circ_1^{[2][n]}(b,x)\large). 
$$
Then $d:H\to H$ is a differential and $H,d\cong D(gp)$.
We thus have a functor $H: \comm-\mbox{modules}\to complexes$
which sends $\O$ to $(H(\O),d)$.
\subsubsection{} The functor $H$ is representable: there exists a 
complex $K$ of $\comm_1$-modules such that $\Hom(K,M)\cong M$
$K$ is freely generated by elements $m_k$, $k\geq 1$ such that
the grading of $m_k$ is $2-k$, each $m_k$ vanish on shuffles,
 and 
$$
dm_{k+1}(x_1,\ldots,x_{k+1})=x_1m_k(x_2,\ldots, x_{k+1})
-(-1)^{k}m_k(x_1x_2,\ldots,x_{k+1})+\ldots
m(x_1,\ldots,x_k)x_{k+1}.
$$

\subsubsection{} Let $\comm_1$ be the operad of commutative algebras with unit.
Via the inclusion
$$
e:\comm\to \comm_1,
 $$
any $\comm_1$-module  $M$ determines a $\comm$-module $e^*M$, and we have a composition $H'=Hg_*$. The advantage of passage 
from $\comm$ to $\comm_1$-modules is that
insertions of $\u$ make the picture more rigid. 
\subsection{} The category of $\comm_1$-modules admits an easy
combinatorial description which we are going to develop.
\subsection{} Note that for any $X$ we have a canonical
element  $m(X)\in\comm_1(X)$  representing  an $(|X|-1)$-fold  product
if $|X|>1$, the identity opeation if $|X|=1$, and the unit if $|X|=0$.
Therefore, for any map $f:X\to Y$ of finite sets we have a canonical
element of $\comm_1(f)$ (see \ref{loadf}) which is
$$
m_f=\bigotimes_{y\in Y}m_{f^{-1}(y)}.
$$ 
let $M$ be a $\comm_1$-module. We can define the operation
$\bigcirc:\comm_1(f)\otimes M(Y)\to M(X)$
by the same rules as in \ref{loadf}.
Note that $m_f\bigcirc m_g\bigcirc x=m_{gf}\bigcirc x$
for any $g:Z\to X$ and $x\in M$. 

Let $f^*:M(Y)\to M(X)$ be such that  $f^*(x)=m_f\bigcirc x$. 
Call this oeration the inverse image.
\subsubsection{Direct image} Let $i:X\to Y$ be an injection of
finite sets.  Let $Z=Y/\Im i$ be the quotient set
and let $z\in Z$ be the image of $X$. We have a canonical identification:
$$
Y\cong X\disjoint Z-\{z\}.
$$
Let $x\in M(X)$. Set 
$$
i_*(x)=\circ_z^{ZX}(m_Z,x)\in  M(Y)\cong M(X\disjoint Z-\{z\}) .
$$
\subsubsection{Partial maps} To describe the inter-relation
between the inverse and direct images introduce the category
$\part$ with finite sets as objects and partially defined maps
as morphisms ( a morphism $f$ from $X$ to $Y$ is a collection
of a subset $U\subset X$ and a map $f|_U:U\to X$, $U$ being called
the domain).
Any injection $i:X\to Y$ uniquely determines  
a quasi-inverse $i^\cdot:Y\to X$ defined on $\Im i\subset Y$
and being inverse on it to $i$. 
Any partial map $f:X\to Y$ with domain $U\subset X$ splits
as $f=f|_{U}i^{\cdot}$, where $i:U\to X$ is the inclusion.
For a $\comm_1$-module $M$ set $f^*=i_*f|_U^*:M(Y)\to M(X)$.
This endows $M$ with a structure of a functor
$\part^{op}\to complexes$. Let $\shpart$ be the category 
of such functors.  One can also say that 
this construction defines a functor from the category
of $\comm_1$-modules to $\shpart$. One sees that this is 
an isomorphism of categories. In the sequel
we do not distinguish between $\comm_1$-modules and
functors from   $\shpart$.

\subsubsection{An analogue of Dold-Puppe theorem} 
The category $\part$ has a lot of projectors. This fact allows us
to simplify the notion of a contravariant functor from $\part$
to $complexes$ in a similar way as the Dold-Puppe theorem
establishes an equivalence between the category of simplicial
vector spaces and the category of complexes.

 Introduce a category $\sur$ as follows.
Objects of $\sur$ are  finite sets and $\morsur(X,Y)$ consists of all surjective
maps $X\to Y$. The empty set is thus an isolated object: it is not connected
by an arrow with  any nonempty
set $X$.

The category $\sur$ is a subcategory of $\part$.
Let $\Shsur$ be the category of  functors  $\sur^{op}\to Complexes$.

We are going to establish an equivalence between $\Shsur$ and 
$\shpart$.
\subsubsection{Notation} set $\ess(X,Y)\subset \M_{\part}(X,Y)$
to be the set of morphisms whose domain is strictly inside $X$.  
\subsubsection{A functor $A:  \shpart\to \Shsur$}
Let $M$ be in $\shpart$. 

Define an object $A(M)\in \Shsur$
by 
$$
A(M)(X)=M(X)/\sum_{i\in\ess(X,Y)}\Im i^*
$$

One sees that the morphisms of subcategory $\sur$  act on 
$A(M)$. Therefore, $A(M)$ is in $\Shsur$.

\subsubsection{The functor $B:\Shsur\to \shpart$}
Let $N$ be a functor $\sur^{op}\to Complexes$. Set 
$$B(N)(X)=\bigoplus\limits_{Y\subset X}N(Y),
$$
where we allow $Y=X$.
for a partial map $f:Z\to X$ and $Y\subset X$ define a map
 $f_Y:N(Y)\to N(f^{-1}Y)$ by setting it 0  if $Y\neq \Im f$ and
$f|_{f^{-1}(Y)}^*$ if $Y=\Im f$. The direct sum of 
$f_Y$ over all subsets of $X$ defines a map
$f^*:B(N(Y))\to B(N(X))$. $B(N)$ is thus an element in $\shpart$.  
\subsubsection{}
\begin{theorem}
 $A$ and $B$ establish an equivalence of the categories
$\shpart$ and $\Shsur$.
\end{theorem}
\subsection{Injective objects in the category of $\comm_1$-modules} 
We are going to develop a homological algebra in the category
of $\comm_1$-modules. We will frequently replace it
with an equivalent category $\Shsur$.
\subsubsection{Forgetful functor} Let $\S$ be the groupoid of finite sets
and $\Sigma$ the category  of functors $\S^{op}\to Complexes$. Note
that $\S^{op}\cong \S$.  The inclusion $\S\to\sur$ induces a functor
$O: \Shsur\to \Sigma$. We want to construct the {\em right} adjoint functor 
$I$ to $O$.  As usual, let $h(X,Y)=k[\morsur(X,Y)]$; $h$ is a functor
$\sur^{op}\times\sur\to \mbox{Vector Spaces}$.  The inclusion
$\S\to \sur$ allows us to consider $h$ as a functor $\S^{op}\times \sur\to
\mbox {Vector Spaces}$.  For ${\frak P}\in\Sigma$ set 
$$
I({\frak P})={\rm Hom}_{\S^{op}}(h,{\frak P}).
$$
We have:
$$
\morsur (M, I({\frak P}))={\rm Mor}_{\S^{op}}\;(M\otimes_\sur h,{\frak P})\cong {\rm Mor}_{\S^{op}}\;(O(M),{\frak P}).
$$
Therefore, $I$ is right adjoint to $O$.
Since $\Sigma$ is semisimple, $I({\frak P})$ are injective for any ${\frak P}$. Since
$M\to IO(M)$ is an inclusion, any object in $\sur$ admits an injective resolution
consisting of objects $I({\frak P})$ for different ${\frak P}$.
\subsubsection{} $I$ commutes with products.
Let ${\frak P}(X)=\prod_i{\frak P}_i(X)$. Then $I({\frak P})(X)\cong \prod_i I({\frak P}_i)(X)$ for any $X$.
\subsubsection{Computation of $B(I({\frak P}))$} We are going to construct the $\comm_1$-module
$BI({\frak P})$ corresponding to $I{\frak P}$. Decompose ${\frak P}=\prod_n {\frak P}_n$, where ${\frak P}_n$ is not zero
only on $n$-element sets.  Consider first the case when ${\frak P}_n([n])=k[S_n]$ 
is a regular $S_n^{op}$-module
Let $V_n=<e_1,\ldots, e_n>$ be the standard $n$-dimensional 
vector space equipped with a $\Z^n$-grading $gr$ such that $gr(e_i)=(0,\ldots,1,\ldots,0)$,  
$1$ being on the $i$-th position. The symmetric group $S_n$ acts naturally on $V$
 
Let $SV$ be the cofree symmetric coalgebra cogenerated by $V$.  Let $\Delta$
be the coproduct and $\epsilon$ the counit. 
Set
$J_n(X)=(\otimes_{x\in X}SV )^{11\cdots 1}$, where the superscript $(11\cdots 1)$
means the corresponding homogeneous part with respect to $gr$.
\subsubsection{Another description of $J_n(X)$} Let $f:[n]\to X$ be a
map. Let 
\begin{equation}\label{sharp}
f^\sharp\in J_n(X)=\otimes_{x\in X} u_x ,
\end{equation}
 where 
$u(x)=\prod_{i\in f^{-1}(x)}e_i\in SV$. It is clear that $ f^{\sharp}$ for different $f$ form 
a basis in $J_n(X)$. 
 
\subsubsection{}
Define the $\comm_1$ structure on $B(I)$ as follows:

-the action of the symmetric group on $J_n(X)$ is by permutations of tensor factors;

-insertion of unit $e_x:J_n(X)\to J_n(X-x)$ is the application of $\epsilon$
to the $x$-th factor;

-Let $D=\{1,2\}$. The composition 
$$\circ_2^{DX}:m\otimes J_n(X)\to J_n(X\disjoint 1)\subset (\otimes_{x\in X}SV )\otimes SV,
$$
is just the multiplication by 1;

-the composition 
$$\circ_x(XD):J_n(X)\otimes m\to J_{n}(X-x\disjoint D)\subset 
\big(\bigotimes_{x'\in X-\{x\}}SV\big)\otimes SV\otimes SV 
$$
is equal to the application of $\Delta$ to the $x$-th factor and putting 
the result into the last two factors.

Thus, $J_n$ is a $\comm_1$-module. The $S_n$ action on $V$ translates into
an $S_n$ action on $J_n$. For an $S_n^{op}$ -module ${\frak P}_n$, define
$J({\frak P}_n)=(J_n\otimes{\frak P}_n)^{S_n}$.  

\begin{Claim} $A(J({\frak P}_n))\cong I({\frak P}_n)$.
\end{Claim}
\pf Note that $A(J({\frak P}_n))(X)\subset J({\frak P}_n)(X)$  is isomorphic via  the map $\sharp$ from (\ref{sharp})
to $(k[\morsur([n],X)]\otimes {\frak P}_n)^{S_n}\cong  {\rm Hom}_{S_n}(\morsur([n],X),{\frak P}_n)\
\cong I({\frak P}_n)(X)$. One can show that this identification is an isomorphism of
functors.
\endpf
\begin{Corollary}  $J({\frak P}_n)\cong BI({\frak P}_n)$.
\end{Corollary}                         
\subsubsection{Computation of $H(J({\frak P}_n))$} 
Note that $H(J({\frak P}_n))\cong (H(J_n)\otimes {\frak P}_n)^{S_n}$.
Furthermore, as a graded vector space, $H(J({\frak P}_n)=FreeLie(SV[-1])^{11\ldots 1}$.
The differential $d$ is compatible with the free Lie algebra structure  and on the 
generator space $SV[-1]$ equals $dx=[1,x]-[\Delta x]$.  The cohomology of such a complex      
is well known
to be the $(1\ldots 1)$-part of the generators space  $V$.   This part is zero
whenever $n\neq 1$ and is one-dimensional when $n=1$.
We have thus proven the following theorem:
\begin{Theorem} $H(J({\frak P}_n))$ is 
acyclic whenever $n\neq 1$;  $H(J({\frak P}_1))$
is quasi-isomorphic to ${\frak P}_1$
\end{Theorem}
\subsubsection{}   For an $M\in \Shsur$, set 
$$
G(M)=\{x\in M([1])| f^*x=0 \mbox{  for any }
f:[n]\to [1]\mbox{ with } n>1.\}
$$
 We have an obvious morphism $i:G(M)\to B(M)([1])\to H(B(M))$
such that $di=0$.  
\begin{Claim} If $M=I({\frak P})$, $i$ is a quasi-isomorphism.
\end{Claim}
\pf Follows directly from the previous section.    
\endpf
\begin{Corollary} 
In the derived category of $\Shsur$, $H$  represents the right  derived functor from $G$.
\end{Corollary}
\subsubsection{Tensor product of $\comm_1$-modules}
We have a natural map $\comm_1\to \comm_1\otimes \comm_1$
which sends $e\in \comm_1(\emptyset)$ to $e\otimes e$ and
$m$ to $m\otimes m$. This map induces a symmetric monoidal structure on the category of $\comm_1$-modules such that the 
tensor product $M\otimes N$ of two modules $M,N$ is defined
by $M\otimes N(X)=M(X)\otimes N(X)$. Below we are describing the 
corresponding symmetric monoidal structure on $\Shsur$  
\subsection{Tensor product in $\Shsur$}  Let $F,G\in Shsur$. Set
$$
F\boxtimes G(X)=\bigoplus\limits_{Y,Z\subset X;X=Y\cup Z}F(Y)\otimes G(Z).
$$ 

For a surjective map $\phi:X'\to X$,$Y,Z\subset X;X=Y\cup Z$,
and $f\otimes g\in F(Y)\otimes G(Z)\subset F\boxtimes G(X)$
set 
$$
\phi^*(f\otimes g)=\sum\limits_{Y',Z'\subset X'; \phi(Y')=Y;\phi(Z')=Z}
\phi|_{Y'}^*(f)\otimes \phi|_{Z'}^*(g).
$$
\subsubsection{}\label{correspondence} We have 
a natural isomorphism $B(F)\otimes B(G)\cong B(F\boxtimes G).$
\subsubsection{Tensor product of $\S$-modules} Let $P,Q$ be functors
from $\S^{op}$ to $Complexes$.  Set $P\boxtimes Q$ be another such a functor described by $P\boxtimes Q(X) =\oplus_{Y\disjoint Z=X}P(Y)\otimes Q(Z)$ with the obvious action of $\S$.
\subsubsection{} We have a  natural isomorphism
$ I(P\boxtimes Q)\to I(P)\boxtimes I(Q)$.
Indeed, 
$$I(P\boxtimes Q)(X)\cong \Hom_{S^{op}}(h_X,P\boxtimes Q)
$$
$$
\cong  \prod \big\{\Hom_{S_n}(\morsur(X,[n]),\bigoplus\limits_{Y\disjoint Z=[n]}
 P(Y)\otimes Q(Z))\big\}
$$
$$
\cong\prod \big\{\Hom_{S_p\times S_q}(X,[p+q]),P(p)\otimes Q(q))
$$
$$
\cong I(P)\boxtimes I(Q).
$$
\subsubsection{} Therefore, we have a natural isomorphism
$$
BI(P\boxtimes Q)\to BI(P)\otimes BI(Q).
$$

\subsubsection{} \label{vanishing}Suppose that $P(\emptyset)=Q(\emptyset)=0$.
Then $P\boxtimes Q([1])=0$. This implies that 
$$
H(BI(P\boxtimes Q))\cong 0,
$$
where $\cong$ means quasi-isomorphic.

If $M,N$ are $\comm_1$ modules with $M(\emptyset)=N(\emptyset)=0$,
one can find  injective resolutions of them consisting only of
the modules $BI(X)$ with $X(\emptyset)=0$. Therefore,
$H(M\otimes N)\cong 0$. Furthermore,  
$$
H(S^kM)=H(\Lambda^k M )=0
$$
 for any $k\geq 2$ , here $S,\Lambda$ stand for symmetric (resp. exterior)
power with respect to $\otimes$.
\subsection{ The $\comm_1$-module $C_\bullet(\g(\cdot))$}
\subsubsection{} The map of operads $\comm_1\to \CC$ defines a 
$\comm_1$-structure on $\CC$ which, in turn, induces    
a $\comm_1$-structure on the linear part $\g(\cdot)$.
The corresponding saction of $\part^{op}$ is defined in 
\ref{direct}-\ref{inverse}.
\subsubsection{} Set $F_k\CC=\C_{\leq k}(\g(\cdot))$. It is a 
$\comm_1$-module filtration.  
\subsubsection{} The associated graded $\comm_1$-module is
isomorphic to $\oplus_n \bigwedge^n\g(\cdot)[n]$. 
Since $\g(\emptyset)=0$ and according to 
\ref{vanishing}, $H(F_k\CC/F_{k-1}\CC) =0$ whenever $k>1$.
\begin{Corollary}
$$
H(\g(\cdot))\to H(\CC)
$$
induces an injection on homology.
\end{Corollary}
To prove the corollary it suffices to note that we have a canonical
splitting $\CC=C_0(\g(\cdot))\oplus C_{>0}(\g(\cdot))$.
Therefore, $(H(\CC))\cong H(C_0)\oplus H(\g(\cdot))$. 
\subsection{Proof of the main theorem} 
To prove our main theorem in the form of \ref{reform}
It suffices to note that $\grt_1\cong H^0(\g(\cdot))$.

\section{Appendix}
\subsection{Operads}\label{operads}
\subsubsection{definition of operad} 
Let $C$ be a symmetric monoidal category. We denote by $\otimes$ the tensor product, by $\u$ the unit object in $C$,  and by $\sigma_{VW}:V\otimes W\to W\otimes V$ the symmetry morphism. 

Let $\S$ be the groupoid of finite sets and their bijections. 
A functor $S\to C$ will be called an $\S$-module.
An {\em operad  without unit} in $C$ is a collection of the  following
data:

-an $\S$-module $\O$;

-for any finite sets $X,Y$ and an element $x\in X$ a morphism 
$$
\circ_x^{X,Y}:\O(X)\otimes \O(Y)\to \O(X-\{x\}\disjoint Y),
$$
natural in $x,X,Y$ ("the insertion onto the $x$-th position"),
such that the following associativity conditions hold:

-For any finite sets $X,Y,Z$ and any $y,z\in X$, $y\neq z$, denote by $\circ_{yz}^{XYZ}$
$$
O(X)\otimes O(Y)\otimes O(Z)\stackrel{\o_y^{XY}\otimes Id}{\to}O(X-\{x\}\disjoint Y)\otimes O(Z)
\stackrel{o_z^{(X-\{x\}\disjoint Y)Z}}\to\O(X-\{y,z\}\disjoint Z).
$$ 
Then $\circ_{yz}^{XYZ}=(Id\otimes \sigma_{ZY})\circ_{zy}^{XZY}$;

-For any finite sets $X,Y,Z$, $y\in X$ and $z\in Y$,
$$
(\circ_y^{X(Y-\{z\}\disjoint Z)}\circ(Id\otimes \circ^{YZ}_z)=\circ_z{(X-\{y\}\disjoint Y)Z}(\circ_y^{XY}\otimes Id).
$$

{\em An operad (with unit)} is an operad without unit endowed with a morphism $\u_e\u\to \O(e)$
for any 1-element set $e$ such that $\u_e$ agree with any change of $e$ by another one-element set $e'$  
and such that the compositions
$$
\O(X)\to \O(X)\otimes \u\stackrel{\Id\otimes \u_e}\to \O(X)\otimes \O(e)\stackrel{\circ_{x}^{Xe}}\to \O(X)
$$
as well as
$$
\O(X)\to \u\otimes \O(X)\stackrel{\u_e\otimes Id}{\to}\O(e)\otimes\O(X)\stackrel{\circ_x^{eX}}\to \O(X)
$$
are the identities.  

 An element $V\in C$ is said
to have {\em an $\O$-algebra structure} if for any finite set

\subsubsection{} \label{loadf}Let $f:X\to Y$ be a map of finite sets.
For an operad $\O$ set 
$$
\O(f)=\bigotimes_{y\in Y} \O(f^{-1}(y)).
$$
Let $Y'\subset Y$ be a subset. 
An element $x\in \O(f)$ is called $Y'$-trivial if $f|_{f^{-1}(Y')}$ is a bijection onto $Y'$
and if
$$
x\in \bigotimes_{y\in Y'} {\rm Im}u_{f^{-1}(y)}\bigotimes_{y\in Y-Y'}\O(f^{-1}(y)).
$$
We have a map $\bigcirc \O(Y)\otimes \O(f)\to \O(X)$ defined by the following conditions:
\begin{enumerate}
\item[-] Let $y\in Y$. Set $Z=Y-y\disjoint f^{-1}(y)$.
Let $F:X\to Z$ be the map such that $F|_{f^{-1}(y)}=Id$; $F|_{X-f^{-1}y}=f$.
For $ x=\bigotimes_{y\in Y}o_y\in \O(f)$ set
$x'=\bigotimes_{y'\in Y-y}o_y'\bigotimes_{x\in f^{-1}(y)} \u$. Then
$\bigcirc(\circ_y^{YZ}(o_y,o),x')=\bigcirc(o,x)$ for any $o\in \O(Y)$.
\item[-] If $f:X\to X$ is identity  and $i=\bigotimes_{x\in X} \u_x$,
then $\bigcirc(f,o)=o$ for any $o\in \O(X)$.
\end{enumerate}
\subsubsection{Shur functor} 
Assume that $C$ have small colimits.
Let $\O$ be an $\S$-module and $V$
an element in $C$. Set $T_V (X)=\bigotimes_{x\in X} V$.
$T_V$ is also an $\S$-module. Set $s(X)=\O(X)\otimes T_V(X)$.
This is again an $\S$-module.  
Set
$$
\Sf_\O(V)=\limdir_{\C}s.
$$
$\Sf$ is a functor $C\to C$ called the {\em Schur functor.}
\subsubsection{Composition}
We claim that $\Sf_{\O_1}\Sf_{\O_2}(V)\cong\Sf_{\A}(V)$, where $\A$
 is defined as follows. Set 
$$
a(X,Y)=\bigoplus_{f:X\to Y} \O_1(Y)\O_2(f)
$$
and $\A(X)=\limdir_\S a(X,\cdot)$.
Set $\O_1\circ \O_2=\A$. One sees that $\circ$ is an associative
tensor product on the category of $\S$-modules.
Note that $\u\circ\O\cong\O\circ\u\cong\O$.
\subsubsection{} If $\O$ is an operad, the map $\bigcirc$ defines
a map $\bigcirc:\O\circ O\to O$.  This makes $\O$ a unital monoid
in the category of $\S$-modules with the tensor structure $\circ$.
\subsubsection{$\O$-algebras} Let $V\in C$ be an element.
$V$ is said to be an $\O$-algebra if
a map $m:\Sf_O(V)\to V$ is prescribed such that $V$ becomes
a unital module over $\O$ viewed as monoid in the category
of $\S$-modules. In other words, the following must be satisfied:
the maps
$$
c_1:\Sf_\O\Sf_\O(V)\to \Sf_{\O\circ O}V\stackrel\bigcirc\to
\Sf_{\O}V\stackrel m\to V
$$
and
$$
c_2:\Sf_\O\Sf_OV\stackrel m\to \Sf_OV\stackrel m\to V
$$
coincide, and
the composition
$$V\cong \Sf_{\u} V\stackrel{\u_e}\to\Sf_{\O}V \stackrel m\to V
$$
is identity.
 \subsubsection{Tensor product of operads} let $\O_1,\O_2$ be operads.
Define an operad $\O_1\otimes \O_2$ by setting
$\O_1\otimes \O_2(X)=\O_1(X)\otimes \O_2(X)$;
$$
(\u_e)_{\O_1\otimes \O_2}=(\u_e)_{\O_1}\otimes(\u_e)_{\O_2};
$$
and 
$$
(\circ^{XY}_{x})_{\O_1\otimes \O_2}=(\circ_x^{XY})_{\O_1}\otimes(\circ_x^{XY})_{\O_2} . 
$$
\subsection{Examples: commutative algebras}
Describe an operad $\comm$ in the category of vector spaces
such that a $\comm$-algebra structure is equivalent 
to a commutative algebra structure.

Set $\comm(X)=k$ to be the constant $\S$-module and
set all the compositions to be the identities.
\subsubsection{} Note that $\comm\otimes \O\cong \O$
for any operad $\O$. In other words, $\comm$ is a
unit in the symmetric monoidal category of operads.
\subsection{Operads \comm\{k\}}.
We are going to define
this operad in the category of graded
vector spaces
such that $\comm\{k\}$-structure on $V$
is equivalent to a structure of commutative algebra on $V[k]$.

\subsubsection{Alternating $S$-module $\Lambda$}
Let $k(X)$ be the vector space spanned by the finite set $X$.
Let $\Lambda(X)=\bigwedge^{|X|}k(X)$. It is clear that
$\Lambda$ is an $S$-module. 

We have natural maps  
$$
\wedge^{XY}:\Lambda(X)\otimes\Lambda(Y)\to
\otimes\Lambda(X\disjoint Y)
$$  and 
$i_x:\Lambda(X)\to\Lambda(X-\{x\})$ for any $x\in X$.
\subsubsection{Operad \comm\{1\}} Set  
$$
\comm\{1\} (X) =\Lambda(X)[|X|-1].
$$ 
Set 
$$
\circ_{x}^{XY}=i_x\circ \wedge^{XY}
$$
and $\u_e={\rm Id}$.
\subsubsection{} Set $\comm\{-1\}$ to be 
an operad such that 
$$
\comm\{\-1\}(X)=\Lambda(X)[1-|X|],
$$
the
composition law being defined by the same formula as for
$\comm\{1\}$
\subsubsection{General case} Set $\comm\{0\}=\comm$;
for $k\neq 0$ set $\comm{k}=\comm{\pm1}^{\otimes\pm k}$, where
the sign is chosen so that $\pm k>0$.
\subsubsection{} For an operad $\O$ set $\O\{k\}$ to be
$\O\otimes \comm\{k\}$. An $\O$-algebra structure on a graded vector space (or complex)
$V$    is equivalent to an $\O$-structure on $V[k]$.
\subsection{Operad of homotopy associative algebras $\hoass$}
\subsubsection{Operad of associative algebras }
Set $\ass_s(X)$ to be the set of all complete orderings on $X$.
Let $X'\in\ass_s(X)$ and  $Y'\in \ass_s(Y)$
For any $x\in X$ introduce a total order on $X-\{x\}\disjoint Y$ by the rule:

-the order on $X-\{x\}$ and $Y$ is the induced one;

-let $\xi\in X-\{x\}$ and $\eta\in Y$. Then $\xi<\eta$ if and only if
$\xi<x$ and $\xi>\eta$ if and only if $\xi>x$.
This rule defines a composition map on $\ass'$. 
Set $\u_e:\u\to \ass_s(e)$ to be the identity. 
This makes $\ass_s$ an operad of sets.
Set $\ass(X)=k[\ass_S(X)]$ to be the span. 
Then $\ass$ is an  operad of vector spaces. An $\ass$-algebra
is the same as an associative algebra. 
\subsubsection{Gerstenhaber bracket in operads} Let $\O$ be an operad.
Let $X,Y,Z$ be finite totally ordered sets such that $\#Z=\#X+\#Y-1$.
For any $x\in X$ introduce a total order on $X-\{x\}\disjoint Y$ 
as was explained in the previous section.
We have then a canonical identification 
$$i_x: X-\{x\}\disjoint Y\to Z.
$$
Let $m_1\in \O(X)$ and $ m_2\in \O(Y)$ set
$m_1\{m_2\}$ to be the following element in $\O(Z)$:
$$
m_1\{m_2\}=\sum\limits_{x\in X} i_{*}\circ_x^{XY}(m_1,m_2).
$$
 Set 
$$
\{m_1,m_2\}=m_1\{m_2\}-(-1)^{|m_1||m_2|}\{m_2,m_1\}.
$$
This operation is called the Gerstenhaber bracket. It satisfy the Jacobi
identity.
\subsubsection{Operad $\hoass'$} Is generated by $n$ ary operations
$m_n'$, $n\geq 2$ such that the grading of each of $m_n'$ is $+1$.
The differential is given by 
$$
dm_n'+\frac 12\sum\limits_{i=2}^{n-1}\{m_i',m_{n+1-i}'\}=0.
$$
In particular, $dm_2'=0$.
\subsubsection{Operad $\hoass$} Set $\hoass=\hoass'{1}$.
denote by $m_n$ the element of  $\hoass$  corresponding to 
 $m_n'\in \hoass'$. The grading of $m_n$ is $2-n$. We have a map
$p:\hoass\to \ass$  such that $p(m_2)=m$ and $p(m_n)=0$ if $n>2$. 
\subsection{Operad of homotopy commutative algebras $\hocomm$}
\subsubsection{Shuffles} \label{shuffles}
let $X=Y\disjoint  Z$ be finite non-empty sets and a total order on $Y$ and $Z$
is prescribed. Set $Sh(Y,Z)$ to be the set of all total orders on $X$
compatible with the ones on $Y,Z$. 
let $sh_{YZ}\in \ass(X)$ be the sum of all elements
from $sh_{YZ}$. Let $sh(X)\in \ass(X)$  be the subspace
generated by all elements of the form $sh_{YZ}$ for all
decompositions $X=Y\disjoint Z$ and all total orders on $Y,Z$.
 For the interval $[n]$, all its total orders are in one-one correspondence
with the elements of the symmetric group $S_n$. 
Let $sh(n)\subset k[S_n]$ be the subspace corresponding to
$sh([n])$. 

We see that $sh\subset \ass$ is an $S$-submodule.
\subsubsection{operad of Lie algebras} Let $\lie$ be the operad of Lie
algebras. It is well known that  the linear dual $S$-module
$\lie^*$ is isomorphic to $\ass/sh$    

\subsubsection{The operad $\hocomm'$} 
 Let $I\in \hoass'$ be the ideal
generated by $\sh(n)m_{n}$.
One checks that $I$ is preserved by the differential.
Set $\hocomm'=\hoass'/I$. 
Note that $\hocomm'$ is a free operad. The $S_n$-submodules
$G'(n)\subset \hocomm'(n)$ generated by $m_n$
are isomorphic to $\lie(n)^*$.
\subsubsection{}\label{hocomm}
Set $\hocomm=\hocomm'{1}$. The map $p$ descends
to a map $p_c:\hocomm\to \comm$, which is  a free resolution of $\comm$.
The images of $m_n$ in $\hocomm$ will be denoted by the same letter.
\subsubsection{} The spaces $G(n)\subset \hocomm(n)$ generated
by the images of $m_n$ are isomorphic to
$\lie\{1\}*$. We will denote $\lie'=\lie\{1\}$.

\end{document}